\documentclass[a4paper,12pt]{article}
\usepackage[francais,english]{babel}
\usepackage{latexsym,amssymb}
\usepackage{amsmath}
\usepackage{pb-diagram}
\newtheorem{thm}{Theorem}[section]

\newtheorem{rema}[thm]{Remark}
\newtheorem{definition}[thm]{Definition}

\def\spin{{\rm{Spin}}}
\def\SO{{\rm{SO}}}

\newcommand{\C}{\mathbb{C}}

\newcommand{\ddim}{\mathrm{dim}}
\newcommand{\lto}{\ensuremath{\longrightarrow}}
\newcommand{\function}[5]
{\begin{eqnarray*}\begin{array}{r@{}ccl}
#1\;\colon\;  & #2 &\lto & #3 \\[.05cm]
  & #4 &\longmapsto  & #5
\end{array}\end{eqnarray*}
}
\begin{document}
\title {{\bf Eigenvalues of the transversal Dirac Operator on 
K{\"a}hler Foliations}}
\author{{\small {\bf {Georges Habib}}}\\
{\small Institut {\'E}lie Cartan,
Universit{\'e} Henri Poincar{\'e}, Nancy I, B.P. 239}\\
{\small 54506 Vand\oe uvre-L{\`e}s-Nancy Cedex, France}\\
{\small {\tt habib@iecn.u-nancy.fr}}}
\date{}
\maketitle
\noindent
\begin{abstract} In this paper,  we
   prove Kirchberg-type inequalities for any K{\"a}hler spin foliation. Their
   limiting-cases are then characterized as being transversal 
minimal Einstein foliations. The key point is to introduce
the transversal K{\"a}hlerian twistor operators. 
\end{abstract}
\section{Introduction} On a compact Riemannian spin manifold $(M^{n},g_{M}),$
Th. Friedrich \cite{F} showed that any eigenvalue $\lambda$ of the Dirac
operator satisfies 
\begin{equation}
\lambda^{2}\geq\frac{n}{4(n-1)}S_{0},
\label{eq:100}
\end{equation}
where $S_{0}$ denotes the infimum of the scalar curvature of $M.$ The limiting case in 
(\ref{eq:100}) is characterized by the existence of a {\it Killing spinor}.
As a consequence $M$ is Einstein. K.D. Kirchberg \cite{K} established that, on such manifolds   
any eigenvalue $\lambda$  satisfies the inequalities
\begin{eqnarray*}
\lambda ^{2}\geq 
\left\{\begin{array}{ll}
\frac{m+1}{4m}S_{0}& \textrm{if $m$ is odd,}\\\\
\frac{m}{4(m-1)}S_{0}& \textrm{if $m$ is even.}
\end{array}\right.
\end{eqnarray*}
On a compact Riemannian spin  foliation $(M,g_{M},\mathcal{F})$ of codimension $q$ with a
bundle-like metric $g_{M}$ such that   the
mean curvature $\kappa$ is a basic coclosed 1-form, S.D. Jung \cite{j1} showed that any eigenvalue $\lambda$ of the transversal Dirac operator satisfies
\begin{equation}
\lambda^{2}\geq\frac{q}{4(q-1)}K_{0}^{\nabla},
\label{eq:569}
\end{equation}      
where $K_{0}^{\nabla}=\mathop{\inf}\limits_M
 (\sigma^{\nabla}+|\kappa|^{2}),$ here  $\sigma^{\nabla}$ denotes the transversal scalar curvature
with the transversal Levi-Civita connection $\nabla.$ The limiting case in (\ref{eq:569}) is
characterized by the fact
that $\mathcal{F}$ is minimal ($\kappa=0$) and transversally
Einstein (see Theorem \ref{thm:002}). 
The main result of this  paper is the following: 
\begin{thm} Let $(M,g_{M},\mathcal{F})$ be a compact Riemannian manifold
  with a K{\"a}hler spin foliation $\mathcal{F}$
  of codimension $q=2m$ and a bundle-like metric $g_{M}.$ Assume that
  $\kappa$ is a basic coclosed 1-form, then  any
  eigenvalue $\lambda$ of the transversal Dirac operator satisfies:
\begin{equation}
\lambda ^{2}\geq \frac{m+1}{4m}K_{0}^{\nabla}\qquad\textrm{if m is odd,}
\label{eq:12}
\end{equation}
and
\begin{equation}
\lambda ^{2}\geq 
\frac{m}{4(m-1)}K_{0}^{\nabla} \qquad\text{if m is even.}\quad
\label{eq:13}
\end{equation}
\end{thm}
The limiting case in \eqref{eq:12}  is characterised  by the fact that the
foliation is minimal and by existence of a transversal K\"ahlerian Killing
spinor (see Theorem \ref{thm:08}). We refer to Theorem \ref{thm:005} for the equality case in \eqref{eq:13}.\\ 
We point out
 that Inequality \eqref{eq:12}  was proved by S. D. Jung  \cite{j2} with
 the additional assumption that
 $\kappa$ is \emph{transversally holomorphic}. 
The author would like to thank Oussama Hijazi for his support. 
\section{Foliated manifolds}
\label{section:2}
\setcounter{equation}{0}
In this section, we  summarize some standard facts about
foliations. For more details, we refer to  \cite{p},
\cite{j1}.\\
Let $(M,g_{M})$ be a $(p+q)$-dimensional Riemannian  manifold
and a foliation $\mathcal{F}$ of codimension $q$ and let $\nabla^{M}$ be
the Levi-civita connection associated with $g_{M}.$  We consider the exact sequence
\begin{equation*}
0  \longrightarrow  L \stackrel{\iota}{\longrightarrow} TM 
\stackrel{\pi}{\longrightarrow} Q \longrightarrow 0 ,
\end{equation*}
where $L$ is the tangent bundle of $TM$ and $Q= TM/L \simeq
L^\perp$ the normal bundle. We assume  $g_{M}$ to be a {\it bundle-like metric} on
$Q$, that means the induced metric $g_{Q}$ verifies the holonomy invariance
condition, 
 \begin{equation*}
 \mathcal{L}_{X}g_{Q}=0,\hskip1cm \forall  X\in  \Gamma(L),
\end{equation*}
where $\mathcal{L}_{X}$ is the Lie derivative with respect to $X$.
Let $\nabla $ be the connection on $Q$ defined by:
\begin{equation*}
\nabla _{X} s =
\left\{\begin{array}{ll}
\pi \left[X,Y_{s}\right ], &  \textrm {$ \qquad{\forall}  X  \in \Gamma(L)$ },\\\\
\pi \left(\nabla_{X}^{M}Y_{s}\right), & \textrm {$ \qquad{\forall} X  \in \Gamma (L^\perp)$ },
\end{array}\right.
\end{equation*}
where $s\in \Gamma(Q)$ and  $Y_{s}$ is the unique vector of $\Gamma(L^\perp)$ 
such that $\pi \left(Y_{s}\right)=s.$
 The connection $\nabla $ is metric and torsion-free. The curvature of
 $\nabla $ acts on $\Gamma(Q)$ by :
$$R^{\nabla} \left(X,Y\right)s=\nabla _{X}\nabla _{Y}s-\nabla _{Y}\nabla 
_{X}s-\nabla _{\left[X,Y\right]}s,\qquad{\forall} X,Y \in \chi 
\left(M\right).$$
The transversal Ricci curvature is defined by:
\function {\rho^{\nabla }}{\Gamma(Q)}{\Gamma(Q)}{X}{\rho^{\nabla 
}\left(X\right)=\displaystyle\sum_{j=1}^{q} R^{\nabla}\left(X,e_{j}\right)e_{j}.}
Also, we define the transversal scalar curvature :
\begin{equation*}
\sigma  ^{\nabla }=\mathop\sum_{i=1}^{q}g_{Q}\left(\rho ^{\nabla 
}\left(e_{i}\right),e_{i}\right)
=\sum_{i,j=1}^{q}R^{\nabla }\left(e_{i},e_{j},e_{j},e_{i}\right),
\end{equation*}
where $ \left\{e_{i}\right\}_{i=1,\cdots, q}$ is a local orthonormal frame
of $Q$ and $R^{\nabla}(X,Y,Z,W)=g_{Q}(R^{\nabla}(X,Y)Z,W),$ for all
$X,Y,Z,W \in \Gamma(Q).$
The foliation  $\mathcal{F}$ is said to be transversally Einstein if and
only if
\begin{equation*}
\rho ^{\nabla }=\frac{1}{q}\sigma ^{\nabla } {\rm Id}, 
\end{equation*}
 with constant transversal scalar curvature.
The mean curvature of $Q$ is given  by:
\begin{equation*}
\kappa \left(X\right)= g_{Q}\left(\tau ,X\right),\hskip2cm \forall X \in 
\Gamma(Q),
\end{equation*}
where $\tau =\sum _{l=1}^{p}II \left(e_{l},e_{l}\right),$
with  $\left \{e_{l}\right\}_{l=1,\cdots, p} $ is a local orthonormal frame
of $\Gamma(L)$ and $II $ is the second fundamental form of  $\mathcal{F}$ 
defined  by:
\function {II}{\Gamma(L) \times \Gamma(L) }{\Gamma(Q)}{\left(X,Y\right)}{II \left(X,Y\right)=\pi \left(\nabla 
^{M}_{X}Y\right).}
We define basic $r$-forms by :
$$\Omega_{B}^{r}\left(\mathcal{F}\right)=\left\{\Phi\in \Lambda ^{r}T^{*}M 
|\hskip0.2cm X\llcorner\Phi =0 \quad\text{and}\quad
X\llcorner d\Phi 
=0,\hskip0.5cm  \forall X \in \Gamma(L) \right\}, $$
where $d$ is the exterior derivative  and 
$X\llcorner$ is the interior product. Any  $\Phi \in
\Omega_{B}^{r}(\mathcal{F})$ can be locally  written as
\begin{equation*}
\sum _{1\leq j_{1}<\cdots<j_{r}\leq q}\beta _{j_{1},\cdots,j_{r}}dy_{j_{1}}\wedge 
\cdots\wedge dy_{j_{r}},
\end{equation*}
where  $\frac{\partial}{\partial x_{l}}\beta_{j_{1},\cdots,j_{r}}=0,\hskip0.3cm  \forall l=1,\cdots, p.$
With the local expression  of  basic $r$-forms, one can verify that
$\kappa$  is closed if $\mathcal{F}$ is isoparametric ($\kappa \in 
\Omega_{B}^{1}\left(\mathcal{F}\right)$). For all  $r\geq 0,$
\begin{equation*}
d\left(\Omega_{B}^{r}\left(\mathcal{F}\right)\right) \subset   
\Omega_{B}^{r+1}\left(\mathcal{F}\right).
\end{equation*}
 We denote by $ d_{B}=d|_{\Omega _{B}\left(\mathcal{F}\right)}$
where $\Omega _{B}\left(\mathcal{F}\right)$ is the tensor algebra of
$\Omega_{B}^{r}\left(\mathcal{F}\right).$
We have the following formulas:
\begin{equation*}
d_{B}=\sum _{i=1}^{q}e_{i}^{\star }\wedge \nabla _{e_{i}} \quad\text{and}\quad  
\delta _{B}=-\sum_{i=1}^{q} {e_{i}}\llcorner \nabla
_{e_{i}}+\kappa\,\llcorner,
\end{equation*}
where  $\delta _{B}$ is the adjoint operator of $ d_{B}$ with respect to the induced
scalar product  and  
$\left\{e_{i}\right\}_{i=1,\cdots,q}$ is a  local orthonormal frame of  $Q$.
\section{The transversal Dirac operator on K\"ahler Foliations}
\setcounter{equation}{0}

In this section, we start by recalling some facts on  Riemannian foliations which
could be found in \cite{gk1}, \cite{gk2}, \cite{ag}, \cite{j1}. For
completeness, we also scketch a straightforward  proof of Inequality (\eqref{eq:569})
established in \cite{j1} and end by recalling well-known facts  (see
\cite{K}, \cite{D}, \cite{h1}, \cite{h2}, \cite{j2}) on  K{\"a}hler spin
foliations.\\

On a foliated Riemannian manifold
  $\left(M,g_{M},\mathcal{F}\right),$ a transversal spin
structure is a pair 
$\left(\spin Q,\eta\right )$ where $\spin Q$ is a  
$\spin_{q}$-principal fibre bundle over $M$ and $\eta $ a $2$-fold cover such that
the following diagram commutes:
$$
\begin{diagram}
\node{\spin Q\times \spin_{q}} \arrow{e}\arrow{s,r}{\eta\otimes Ad}
\node{\spin Q} \arrow{s,r}{\eta} \arrow{e,t}{} \node{M} \\
\node{\SO Q\times \SO_q}\arrow{e} \node{\SO Q} \arrow{ne,b}{}
\end{diagram}
$$
The maps
$\spin Q\times \spin_{q}\longrightarrow \spin Q,$
and
$\SO Q\times SO_q\longrightarrow \SO Q,$
are respectively the  actions of  $\spin_q$ and $\SO_q$ on the principal
fibre bundles $\spin Q$ and $\SO Q$.
In this case,  $\mathcal{F}$ is called a  transversal spin foliation.
We define the  foliated spinor bundle by:
$ S\left(\mathcal{F}\right):=\spin Q\times _{\rho }\Sigma _{q},$
where
$\rho :\spin_{q}\longrightarrow {\rm Aut}\left(\Sigma_{q}\right),$
is the complex spin  representation and  $
\Sigma_{q}$ is a $\C$ vector space of dimension $N$ with
   $N=2^{\left[\frac{q}{2}\right]},$ where $[\;]$ stands for
the integer part. Recall that the  Clifford  multiplication $\mathcal{M}$ 
on $S\left(\mathcal{F}\right)$ is
 given  by:\function{\mathcal{M}}
{\Gamma(Q)\times \Gamma(S(\mathcal{F}))} {\Gamma(S(\mathcal{F}))}
{(X,\Psi)}{X\cdot \Psi.}

There is a natural Hermitian product on $S\left(\mathcal{F}\right)$ such
that, for all  $X,Y \in \Gamma(Q)$, the following relations are true:
\begin{eqnarray*}
\left<X\cdot \Psi ,\Phi\right >&=&-\left<\Psi ,X\cdot \Phi\right >,\\ 
X\left(\left<\Psi ,\Phi\right >\right)&=&\left<\nabla _{X}\Psi ,\Phi\right 
 >+\left<\Psi ,\nabla _{X}\Phi\right >,\\
\nabla _{Y}\left(X\cdot \Psi\right)& =&\left(\nabla _{Y}X\right)\cdot \Psi  
+X\cdot\left(\nabla _{Y}\Psi\right),
\end{eqnarray*}
where $\nabla$ is the Levi-Civita connection on
$S\left(\mathcal{F}\right)$ and $\Psi,\Phi \in\Gamma( S(\mathcal{F})).$ \\ 

The transversal Dirac operator \cite{gk1,gk2} is locally given by:
\begin{equation}
D_{tr}\Psi =\sum_{i=1}^{q}e_{i}\cdot\nabla  _{e_{i}}\Psi -\frac{1}{2} 
\kappa \cdot\Psi,
\label{eq:150}
\end{equation}
 for all $\Psi \in \Gamma(S(\mathcal{F})).$ We can easily prove using Green's
 theorem \cite{yt} that this operator is formally self
 adjoint. Furthermore, in \cite{gk2} it is proved that if  $\mathcal{F}$ 
is isoparametric and $\delta_{B} 
\kappa =0$, then we have the Schr{\"o}dinger-Lichnerowicz formula:
$$D_{tr}^{2}\Psi =\nabla _{tr}^{\star }\nabla _{tr}\Psi 
+\frac{1}{4}K_{\sigma}^{\nabla   }\Psi, $$
where $ K_{\sigma  }^{\nabla   }=\sigma^{\nabla}+|\kappa|^{2} $ and
$$ \nabla _{tr}^{\star }\nabla _{tr}\Psi =-\sum_{i=1}^{q}\nabla _{e_{i},e_{i}}^{2}\Psi 
+\nabla _{\kappa}\Psi, $$ with $
\nabla_{X,Y}^{2}=\nabla_{X}\nabla_{Y}-\nabla_{\nabla_{X}Y},$ for all $X,Y
\in \Gamma(TM).$ 
Denote by  $\mathcal{P}$ the transversal twistor operator defined by
\begin{equation*}\mathcal{P}: \Gamma(S(\mathcal{F}))\stackrel{\nabla^{tr}}\longrightarrow \Gamma(Q^{*}\otimes
S(\mathcal{F}))\stackrel{\pi}\longrightarrow\Gamma(\ker \mathcal{M}),
\end{equation*}
where $\pi$ is the orthogonal  projection on the kernel of the Clifford
multiplication $\mathcal{M}$. With respect to a local orthonormal frame $\{e_{1},\cdots,e_{q}\}$, for all
$\Psi \in \Gamma(S(\mathcal{F}))$, one has 
\begin{equation} 
\mathcal{P}\Psi=\sum_{i=1}^{q}e_{i}^{*}\otimes(\nabla_{e_{i}}\Psi+\frac{1}{q}e_{i}\cdot
D_{tr}\Psi+\frac{1}{2q} e_{i}\cdot \kappa\cdot\Psi).
\label{eq:25}
\end{equation}

For any spinor field $\Psi,$ one can easily show that
\begin{equation}
\sum_{i=1}^{q}e_{i}\cdot \mathcal{P}_{e_{i}}\Psi=0.
\label{eq:654}
\end{equation} 
Now we give a simple proof of the following theorem:
\begin{thm} \label{thm:002}\emph{\cite{j1}} 
 Let $\left( M,g_{M},\mathcal{F}\right)$ be a compact Riemannian manifold
 with a spin foliation  $\mathcal{F}$ of  codimension $q$ and a bundle-like
 metric $g_{M}$ with $\kappa \in\Omega _{B}^{1}(\mathcal{F}).$ Assume that
 $\delta_{B}\kappa=0$ and let  $\lambda$ be an eigenvalue of the transversal
 Dirac operator, then 
\begin{equation}
\lambda ^{2}\geq \frac{q}{4(q-1)}K_{0}^{\nabla}.
\label{eq:43}
\end{equation} 
\end{thm}
{\bf Proof.} For all $\Psi \in \Gamma(S(\mathcal{F})),$ we have using
Identities \eqref{eq:25}, \eqref{eq:654} \eqref{eq:150}, 
\begin{equation*}
|\mathcal{P}\Psi|^{2}=|\nabla^{tr}\Psi|^{2}-\frac{1}{q}|D_{tr}\Psi|^{2}-\frac{1}{q}\Re(D_{tr}\Psi,\kappa\cdot\Psi)-\frac{1}{4q}|\kappa|^{2}|\Psi|^{2}.
\end{equation*}
For any spinor field $\Phi,$ we have that
$(\Phi,\kappa\cdot\Phi)=-(\kappa\cdot\Phi,\Phi)=-\overline{(\Phi,\kappa\cdot\Phi)},$
so  the  scalar product $(\Phi,\kappa\cdot\Phi)$ is a pure imaginary function.
Hence for any eigenspinor $\Psi$ of the transversal Dirac operator, we obtain
\begin{equation*}
\int_{M}|\mathcal{P}\Psi|^{2}+\frac{1}{4q}\int_{M}|\kappa|^{2}|\Psi|^{2}=\int_{M}|\nabla^{tr}\Psi|^{2}-\frac{1}{q}\int_{M}\lambda^{2}|\Psi|^{2},
\end{equation*}
from which we deduce \eqref{eq:43}  with the help of  the Schr{\"o}dinger-Lichnerowicz formula. Finally, we
can easily prove in the limiting case that $\mathcal{F}$ is minimal
i.e. $\kappa=0,$  and transversally Einstein.
\hfill$\square$ \\
 
A foliation  $\mathcal{F}$ is called K{\"a}hler 
  if there exists a complex parallel orthogonal structure $J:\Gamma(Q)\longrightarrow
  \Gamma(Q)\;(\ddim Q=q=2m)$.
Let $\Omega$ be the associated K{\"a}hler, i.e., for all $X,Y \in \Gamma(Q),$ 
$\Omega(X,Y)=g_{Q}(J(X),Y)=-g_{Q}(X,J(Y))$.
The K{\"a}hler form can be locally expressed  as 
\begin{equation*}
\Omega = \frac{1}{2}\sum_{i=1}^{q}e_{i}\cdot J(e_{i})=
-\frac{1}{2}\sum_{i=1}^{q}J(e_{i})\cdot e_{i},
\end{equation*}
and for all $X \in \Gamma(Q)$, we have $[\Omega,X]:=\Omega\cdot X-X\cdot\Omega=2J(X).$
Under the action of the K{\"a}hler form, the
  spinor bundle splits into an orthogonal sum  
\begin{equation*}
 S(\mathcal{F})=\mathop\oplus_{r=o}^{m} S_{r}(\mathcal{F}),
\end{equation*}
where $ S_{r}(\mathcal{F})$ is an  eigenbundle associated with the eigenvalue
$i\mu_{r}=i(2r-m)$ of the K{\"a}hler form $\Omega.$
Moreover, the spinor bundle of a K{\"a}hler spin foliation carries a
parallel  anti-linear
map $j$ satisfying the relations:
\begin{eqnarray*}
j^{2} &=& (-1)^\frac{m(m+1)}{2} Id,\\
\small[X,j\small] &=& 0,\\
(j\Psi,j\Phi)&=& (\Phi,\Psi),
\end{eqnarray*}
and we have $j\Psi_{r}=(j\Psi)_{m-r}.$
For all $X \in \Gamma(Q),$ we have 
\begin{equation*}
  p_{+}(X)\cdot S_{r}(\mathcal{F})\subset
  S_{r+1}(\mathcal{F})\quad\text{and}\quad
  p_{-}(X)\cdot S_{r}(\mathcal{F})\subset S_{r-1}(\mathcal{F}),
\end{equation*}
 where $
  p_{\pm}(X)=\frac{X \mp iJ(X)}{2}.$
 We define the operator $\widetilde {D}_{tr}$ by
\begin{equation*} 
\widetilde {D}_{tr}\Psi =\sum_{i=1}^{q}J(e_{i})\cdot\nabla _{e_{i}}\Psi
 -\frac{1}{2}J(\kappa)\cdot\Psi.
\end{equation*}
 The local expression of $\widetilde {D}_{tr}$  is independant of the choice
 of the local frame and by Green's theorem \cite{yt}, we
 prove that this operator is self-adjoint.
 On a K{\"a}hler spin foliation, the operators $D_{tr}$ and $\widetilde{D}_{tr}$ satisfy:
\begin{eqnarray}
\small[\Omega,D_{tr}\small]&=&2\widetilde{D}_{tr},\\
\small[\Omega,\widetilde{D}_{tr}\small]&=&-2D_{tr},\\
\small[\Omega,D_{tr}^{2}\small]&=&0,\label{eq:698}\\
D_{tr}\widetilde{D}_{tr}+\widetilde{D}_{tr}D_{tr}&=&0,\label{eq:699}\\
\widetilde{D}_{tr}^{2}&=&D_{tr}^{2}\label{eq:700}.
\end{eqnarray}
We should point out that Equations \eqref{eq:698}, \eqref{eq:699} and
\eqref{eq:700} are true under the assumptions that $\mathcal{F}$ is
isoparametric and $\delta_{B}\kappa=0.$  Now we define the two operators $D_{+}$ and $D_{-}$ by 
\begin{equation}
D_{+}=\frac{1}{2}(D_{tr}-i\widetilde{D}_{tr}) \quad\text{and}\quad
D_{-}=\frac{1}{2}(D_{tr}+i\widetilde{D}_{tr}).
\label{eq:1894}
\end{equation}
 Furthermore, $D_{tr}$ splits into $D_{+}$ and $D_{-},$ and we have the two exact
sequences:
\begin{eqnarray}
\Gamma(S_{m}(\mathcal{F}))\stackrel{D_{-}}\longrightarrow \ldots
\Gamma(S_{r}(\mathcal{F}))\stackrel{D_{-}}\longrightarrow\Gamma(S_{r-1}(\mathcal{F}))\stackrel{D_{-}}\longrightarrow
\ldots \Gamma(S_{0}(\mathcal{F})),\label{eq:1895}\\\nonumber\\
\Gamma(S_{0}(\mathcal{F}))\stackrel{D_{+}}\longrightarrow \ldots
\Gamma(S_{r}(\mathcal{F}))\stackrel{D_{+}}\longrightarrow\Gamma(S_{r+1}(\mathcal{F}))\stackrel{D_{+}}\longrightarrow
\ldots \Gamma(S_{m}(\mathcal{F})).\label{eq:1896}
\end{eqnarray}
\section{Eigenvalues of the transversal Dirac operator}
 \setcounter{equation}{0}
In this section, we prove Kirchberg-type inequalities  by using the
transversal  K{\"a}hlerian
 twistor operators on  K{\"a}hler spin foliations. We refer to \cite{k2}, \cite{k1}.
\begin{definition}On a K{\"a}hler spin foliation, we define the transversal  K{\"a}hlerian
  twistor operators by 
\begin{equation*}
 \mathcal{P}^{(r)}:\Gamma(S_{r}(\mathcal{F}))\stackrel{\nabla^{tr}}\longrightarrow
\Gamma(Q^{*}\otimes S_{r}(\mathcal{F}))\stackrel{\pi_{r}}\longrightarrow
\Gamma(\ker \mathcal{M}_{r}),
\end{equation*}
where $\mathcal{M}_{r}$ is the transversal Clifford multiplication defined by
 \function{\mathcal{M}_{r}}{\Gamma (Q^{*}\otimes S_{r}(\mathcal{F}))}{\Gamma(
  S_{r-1}(\mathcal{F}))\oplus\Gamma
  (S_{r+1}(\mathcal{F}))}{X\otimes\Psi_{r}}{p_{-}(X)\cdot\Psi_{r}\oplus
  p_{+}(X)\cdot\Psi_{r}.}
\end{definition}
For all $r \in \{0,\ldots,m\}$ and  $\Psi_{r} \in
\Gamma(S_{r}(\mathcal{F}))$, we have 
\begin{equation}
\mathcal{P}^{(r)}\Psi_{r}=\sum_{i=1}^{q}
e_{i}^{*}\otimes(\nabla_{e_{i}}\Psi_{r}+a_{r} p_{-}(e_{i})\cdot
\mathcal{D}_{+}\Psi_{r}+b_{r} p_{+}(e_{i})\cdot \mathcal{D}_{-}\Psi_{r}),
\label{eq:111}
\end{equation}
where $\mathcal{D}_{\pm}=D_{\pm}+ \frac{1}{2} p_{\pm}(\kappa) \quad\text{with}\quad a_{r}=\frac{1}{2(r+1)} \quad\text{and}\quad b_{r}=\frac{1}{2(m-r+1)}.$
For any spinor field $\Psi_{r}\in \Gamma(S_{r}(\mathcal{F})),$ we can easily prove 
\begin{equation}
 \sum_{i=1}^q e_{i}\cdot\mathcal{P}^{(r)}_{e_{i}}\Psi_{r}=0.
\label{eq:120}
\end{equation} 
\begin{rema} For any non zero eigenvalue $\lambda$ of
  $D_{tr}$, there exists a spinor field $\Psi \in \Gamma(S(\mathcal
  F))$ called of type $(r,r+1),$ such that 
$D_{tr}\Psi=\lambda \Psi \quad\text{and}\quad \Psi=\Psi_{r}+\Psi_{r+1},$ 
with $r\in \{0,\cdots ,m-1\}$. 
By using \eqref{eq:1894}, \eqref{eq:1895} and
  \eqref{eq:1896} it follows that
  $D_{-}\Psi_{r}=D_{+}\Psi_{r+1}=0,$
  $D_{-}\Psi_{r+1}=\lambda\Psi_{r},D_{+}\Psi_{r}=\lambda\Psi_{r+1}$ and $\|\Psi_{r}\|_{L^{2}}=\|\Psi_{r+1}\|_{L^{2}}.$  
\end{rema}
{\bf Proof.} Let $\varphi$ be an eigenspinor of $D_{tr}.$  There exists an
$r$ such that  $\varphi_{r}$ does not vanish. Let
$\Psi=\frac{1}{\lambda}D_{-}D_{+}\varphi_{r}+D_{+}\varphi_{r},$ one can
easily get that $D_{tr}\Psi=\lambda\Psi.$
\begin{thm}\label{thm:08}  Let $\left( M,g_{M},\mathcal{F}\right)$ be a  compact Riemannian
manifold with a K{\"a}hler spin foliation $\mathcal{F}$ of codimension $q=2m$
and a bundle-like metric
$g_{M}$ with $\kappa\in \Omega_{B}^{1}(\mathcal{F})$ and
$\delta_{B}\kappa=0.$ Then any eigenvalue $\lambda$ of the transversal Dirac
operator, satisfies 
\begin{equation}
\lambda ^{2}\geq 
\frac{m+1}{4m}K_{0}^{\nabla}. 
\label{eq:02}
\end{equation}
If $\Psi$ is an eigenspinor of type $(r,r+1)$ associated with an eigenvalue
$\lambda$ satisfying equality in \eqref{eq:02}, then  $r=\frac{m-1}{2},$ 
the foliation $\mathcal{F}$ is minimal and  for all $X\in \Gamma(Q),$ the
spinor $ \Psi$ 
 satisfies 
\begin{equation}
\nabla_{X}\Psi +\frac{\lambda}{2(m+1)}(X\cdot
\Psi-i\varepsilon J(X)\cdot\bar{\Psi})=0,
\label{eq:06}
\end{equation}
 where $\varepsilon=(-1)^{\frac{m-1}{2}},$ and $\bar\Psi:=(-1)^{r}(\Psi_{r}-\Psi_{r+1}).$ As a consequence $m$ is odd and $\mathcal{F}$ is transversally Einstein
with non negative constant transversal curvature $\sigma^{\nabla}.$
\end{thm}
{\bf Proof.} For all $\Psi_{r}\in \Gamma(S_{r}(\mathcal F)),$ 
using Identities \eqref{eq:111} and \eqref{eq:120}, we have
\begin{eqnarray*}\displaystyle
|\mathcal{P}^{(r)}\Psi_{r}|^{2}&=&\sum_{i=1}^{q}|\mathcal{P}^{(r)}_{e_{i}}\Psi_{r}|^{2}=\sum_{i=1}^{q}(\mathcal{P}^{(r)}_{e_{i}}\Psi_{r},\nabla_{e_{i}}\Psi_{r})\\
&=&\sum_{i=1}^{q}(\nabla_{e_{i}}\Psi_{r}+ a_{r}p_{-}(e_{i})\cdot
\mathcal{D}_{+}\Psi_{r}\\&&+b_{r}p_{+}(e_{i})\cdot
\mathcal{D}_{-}\Psi_{r},\nabla_{e_{i}}\Psi_{r}).\\
\end{eqnarray*}
Finally we obtain,
\begin{equation}
|\mathcal{P}^{(r)}\Psi_{r}|^{2}=|\nabla^{tr}\Psi_{r}|^{2}-a_{r}|\mathcal{D}_{+}\Psi_{r}|^{2}-b_{r}|\mathcal{D}_{-}\Psi_{r}|^{2}.
\label{eq:60}
\end{equation}
 Let $\lambda$ be an eigenvalue of $D_{tr}$ and let $\Psi$ an eigenspinor
 of type $(r,r+1).$ Applying  
 Equality \eqref{eq:60}  to $\Psi_{r},$ one gets 
\begin{eqnarray*}
|\mathcal{P}^{(r)}\Psi_{r}|^{2}&=&|\nabla^{tr}\Psi_{r}|^{2}-a_{r}\lambda^{2}|\Psi_{r+1}|^{2}-a_{r}\lambda
\Re(\Psi_{r+1},p_{+}(\kappa)\cdot\Psi_{r})\nonumber\\&&
-\frac{a_{r}}{4}|p_{+}(\kappa)\cdot\Psi_{r}|^{2}-\frac{b_{r}}{4}|p_{-}(\kappa)\cdot\Psi_{r}|^{2}.
\end{eqnarray*}
By the Schr{\"o}dinger-Lichnerowicz formula and by the fact that $\Psi_{r}$
and $\Psi_{r+1}$ have the same $L^{2}$-norms, we get 
\begin{eqnarray}
\int_{M}|\mathcal{P}^{(r)}\Psi_{r}|^{2}+\frac{a_{r}}{4}\int_{M}|p_{+}(\kappa)\cdot\Psi_{r}|^{2}+\frac{b_{r}}{4}\int_{M}|p_{-}(\kappa)\cdot\Psi_{r}|^{2}\nonumber
&=&\\
\int_{M}((1-a_{r})\lambda^{2}-\frac{1}{4}K_{\sigma}^{\nabla})|\Psi_{r}|^{2}-a_{r}\lambda
\int_{M}\Re(\Psi_{r+1},p_{+}(\kappa)\cdot\Psi_{r}).
\label{eq:4567}
\end{eqnarray}
Similarly applying  \eqref{eq:60} to $\Psi_{r+1},$ we obtain
\begin{eqnarray}
\int_{M}|\mathcal{P}^{(r+1)}\Psi_{r+1}|^{2}+\frac{a_{r+1}}{4}\int_{M}|p_{+}(\kappa)\cdot\Psi_{r+1}|^{2}+\frac{b_{r+1}}{4}\int_{M}|p_{-}(\kappa)\cdot\Psi_{r+1}|^{2}\nonumber
&=&\\
\int_{M}((1-b_{r+1})\lambda^2-\frac{1}{4}K_{\sigma}^{\nabla})|\Psi_{r+1}|^{2}+b_{r+1}\lambda
\int_{M}\Re(\Psi_{r+1},p_{+}(\kappa)\cdot\Psi_{r}), 
\label{eq:4568}
\end{eqnarray} 
where $K_{\sigma}^{\nabla}=\sigma^{\nabla}+|\kappa|^{2}.$ In order to get
rid  the term $\lambda
\int_{M}\Re(\Psi_{r+1},p_{+}(\kappa)\cdot\Psi_{r}),$  since the
l.h.s. of \eqref{eq:4567} and \eqref{eq:4568} are non negative,  dividing
 \eqref{eq:4567} by $ a_{r}$ and  \eqref{eq:4568} by $b_{r+1}$ then
 summing up,  we find  by substituting the values of  $a_{r}$ and $b_{r+1},$
\begin{equation*}
\lambda^{2}\geq \frac{m+1}{4m}K_{0}^{\nabla}.
\end{equation*}
Now, we discuss the limiting case of Inequality \eqref{eq:02}. Dividing  \eqref{eq:4567} by $ a_{r}$ and  \eqref{eq:4568} by $b_{r+1}$ then
 summing up as before, and substituting  $a_{r},$ $b_{r+1}$ and $\lambda^{2}$ by
 their values, we easily deduce that $\kappa=0,$
 $\mathcal{P}^{(r)}\Psi_{r}=0$ and
 $\mathcal{P}^{(r+1)}\Psi_{r+1}=0.$ Hence by  \eqref{eq:4567}, we
 find that
$
\lambda^{2}=\frac{1}{4(1-a_{r})}\sigma_{0}=\frac{m+1}{4m}\sigma_{0}$ where
$\sigma_{0}=\mathop{\inf}\limits_{M}\sigma^{\nabla},$ then $r=\frac{m-1}{2}$ and $m$ is odd.
 It remains to prove that $\Psi$ satisfies \eqref{eq:06}. For
 $r=\frac{m-1}{2},$  by definition of the K{\"a}hlerian twistor operators, for all $j\in\{1,\cdots,q\},$ we obtain
\begin{equation*}
\nabla_{e_{j}}\Psi_{r}+\frac{\lambda}{m+1}p_{-}(e_{j})\cdot\Psi_{r+1}=0,
\end{equation*}
and
\begin{equation*}
\nabla_{e_{j}}\Psi_{r+1}+\frac{\lambda}{m+1}p_{+}(e_{j})\cdot\Psi_{r}=0.
\end{equation*}
 Summing up the two equations, we get
\eqref{eq:06} for $X=e_{j}.$  Using  Ricci identity in \eqref{eq:06}, one  easily proves that $\mathcal{F}$ is transversally
Einstein. 
\hfill$\square$  
\begin{thm} \label{thm:005} Under the same conditions as in Theorem \emph{\ref{thm:08}}
 for $m$  even,  any eigenvalue $\lambda$ of the transversal Dirac
  operator satisfies
\begin{equation}
\lambda^{2}\geq \frac{m}{4(m-1)}K_{0}^{\nabla}.
\label{eq:554}
\end{equation}
If $\Psi$ is an eigenspinor of type $(r,r+1)$ associated with an eigenvalue
satisfying equality in \eqref{eq:554}, then   $r=\frac{m}{2},$ the foliation
$\mathcal{F}$ is minimal and $\Psi$ satisfies for all $X \in \Gamma(Q),$
\begin{equation} 
\nabla_{X}\Psi_{r+1}=-\frac{\lambda}{q}(X-iJX)\cdot \Psi_{r}.
\label{eq:008}
\end{equation}
\end{thm}
{\bf Proof.} Let $\Psi$ an eigenspinor of
type $(r,r+1)$ associated with any eigenvalue $\lambda$ of the transversal
Dirac operator $D_{tr}.$ Recalling Equalities \eqref{eq:4567} and \eqref{eq:4568}, we
have
\begin{equation}
 0\leq \int_{M}((1-a_{r})\lambda^{2}-\frac{1}{4}K_{\sigma}^{\nabla})|\Psi_{r}|^{2}-a_{r}\lambda
\int_{M}\Re(\Psi_{r+1},p_{+}(\kappa)\cdot\Psi_{r}),
\label{eq:555}
\end{equation}
and
\begin{equation}
 0\leq \int_{M}((1-b_{r+1})\lambda^{2}-\frac{1}{4}K_{\sigma}^{\nabla})|\Psi_{r+1}|^{2}+b_{r+1}\lambda
\int_{M}\Re(\Psi_{r+1},p_{+}(\kappa)\cdot\Psi_{r}).
\label{eq:556}
\end{equation}
Hence if $\lambda
\int_{M}\Re(\Psi_{r+1},p_{+}(\kappa)\cdot\Psi_{r})\leq 0$, then by \eqref{eq:556} 
\begin{equation*}
\lambda^{2}\geq \frac{1}{4(1-b_{r+1})}K_{0}^{\nabla},
\end{equation*}
The antilinear isomorphism $j$ sends $S_{r}(\mathcal{F})$ to $S_{m-r}(\mathcal{F}).$ This
allows the choice of $\mu_{r}$ to be non negative (i.e. $r\geq \frac{m}{2}$) where $\mu_{r}$ is the eigenvalue
associated with $\Psi_{r}.$  Then a careful study of the graph of the function
$\frac{1}{1-b_{r+1}},$ yields \eqref{eq:554}. \\\\
On the other
hand if $\lambda
\int_{M}\Re(\Psi_{r+1},p_{+}(\kappa)\cdot\Psi_{r})> 0.$  Applying
 Equation (\ref{eq:60}) to the spinor $j\Psi,$ which is a spinor of type
 $(m-(r+1),m-r),$ we find the same inequalities as
 \eqref{eq:555} and \eqref{eq:556}, then 
\begin{equation*}
\lambda^{2}> \frac{1}{1-a_{r}}\frac{K_{0}^{\nabla}}{4}.
\end{equation*} 
  As before we can choose
$\mu_{m-(r+1)}\geq 0$ (i.e. $r \leq \frac{m}{2}-1$). A
careful study of the graph of the function $\frac{1}{1-a_{r}}$ gives
Inequality \eqref{eq:554}.\\
Now we discuss the limiting case of \eqref{eq:554}. As we have seen, it could not be achieved if $\lambda
\int_{M}\Re(\Psi_{r+1},p_{+}(\kappa)\cdot\Psi_{r})> 0,$ so only
the other case should be considered. By \eqref{eq:4568}, one has 
$$
\begin{array}{llc}
\int_{M}|\mathcal{P}^{(r+1)}\Psi_{r+1}|^{2}+\frac{a_{r+1}}{4}\int_{M}|p_{+}(\kappa)\cdot\Psi_{r+1}|^{2}
\\\\+\frac{b_{r+1}}{4}\int_{M}|p_{-}(\kappa)\cdot\Psi_{r+1}|^{2}-b_{r+1}\lambda
\int_{M}\Re(\Psi_{r+1},p_{+}(\kappa)\cdot\Psi_{r})\nonumber=\\\\
(1-b_{r+1})\int_{M}(\frac{m}{4(m-1)}K_{0}^{\nabla}-\frac{1}{4(1-b_{r+1})}K_{\sigma}^{\nabla})|\Psi_{r+1}|^2. 
\label{eq:002}
\end{array}
 $$
Since
$\frac{m}{m-1}=\mathop{\inf}\limits_{r\geq
  \frac{m}{2}}\frac{1}{1-b_{r+1}},$  and the l.h.s. of \eqref{eq:002} is
non negative, we deduce that $\kappa=0, \mathcal{P}^{r+1}\Psi_{r+1}=0$
and  $\frac{m}{m-1}=\frac{1}{1-b_{r+1}}$ so  $r=\frac{m}{2}.$   It remains
to show that  Equation \eqref{eq:008} holds. For this, take $X=e_{j}$ where $\{e_{j}\}_{j=1,\cdots,q}$ is a local
orthonormal frame. For $r=\frac{m}{2},$  and  by definition of the
K{\"a}hlerian twistor operators, for all $j \in \{1,\cdots, q\}$, we obtain
\begin{equation*} 
\nabla_{e_{j}}\Psi_{r+1}+\frac{\lambda}{q}(e_{j}-iJe_{j})\cdot\Psi_{r}=0.
\end{equation*}    
\hfill$\square$
\bibliographystyle{amsalpha}
\bibliography{Article1}
\end{document}